\documentclass{amsart}
\usepackage{amsmath, amssymb, amsthm,amsfonts}
\usepackage{graphicx,color}
\usepackage{graphics}
\usepackage{comment}
\usepackage[OT2,T1]{fontenc}
\DeclareSymbolFont{cyrletters}{OT2}{wncyr}{m}{n}
\DeclareMathSymbol{\Sha}{\mathalpha}{cyrletters}{"58}

\DeclareMathOperator{\SL}{SL}
\DeclareMathOperator{\GL}{GL}
\DeclareMathOperator{\Sh}{Sh}

\DeclareMathOperator{\End}{End}

\newcommand{\Hup}{\mathbb{H}}
\newcommand{\TT}{\mathbb{T}}

\newcommand{\p}{{\mathfrak P}}

\newcommand{\Q}{{\mathbb Q}}
\newcommand{\Z}{{\mathbb Z}}

\newcommand{\C}{{\mathbb C}}
\newcommand{\R}{{\mathbb R}}

\newcommand{\kro}[2]{\left( \frac{#1}{#2} \right) }

            \DeclareFontFamily{U}{wncy}{} 
            \DeclareFontShape{U}{wncy}{m}{n}{% 
               <5>wncyr5% 
               <6>wncyr6% 
               <7>wncyr7% 
               <8>wncyr8% 
               <9>wncyr9% 
               <10>wncyr10% 
               <11>wncyr10% 
               <12>wncyr6% 
               <14>wncyr7% 
               <17>wncyr8% 
               <20>wncyr10% 
               <25>wncyr10}{} 
\DeclareMathAlphabet{\cyr}{U}{wncy}{m}{n}

\begin {document}

\newtheorem{thm}{Theorem}

\newtheorem{lem}{Lemma}[section]
\newtheorem{prop}[lem]{Proposition}

\newtheorem{cor}[lem]{Corollary}

\theoremstyle{definition}

\theoremstyle{remark}

\title[Hecke operators in half-integral weight]{
Hecke operators in half-integral weight
}

\author{Soma Purkait}
\address{Mathematics Institute\\
	University of Warwick\\
	Coventry\\
	CV4 7AL \\
	United Kingdom}

\email{S.Purkait@warwick.ac.uk}
\date{\today}
\thanks{The author is supported by a Warwick Postgraduate Research Scholarship
and an EPSRC Delivering Impact Award}


\begin{abstract}
In \cite{Shimura}, Shimura introduced modular forms of half-integral weight, their Hecke algebras
and their relation to integral weight modular forms via the Shimura correspondence. For modular 
forms of integral weight, Sturm's bounds give generators of the Hecke algebra as a module.
We also have well-known recursion formulae for the operators $T_{p^\ell}$ with $p$ prime. 
It is the purpose of this paper to prove analogous results in the half-integral weight
setting. We also give an explicit formula for how operators $T_{p^{\ell}}$ commute
with the Shimura correspondence.
\end{abstract}
\maketitle

\section{Introduction}

In \cite{Shimura}, Shimura introduced modular forms of half-integral weight, their Hecke algebras
and their relation to integral weight modular forms via the Shimura correspondence. For modular 
forms of integral weight, Sturm's bounds give generators of the Hecke algebra as a module.
We also have well-known recursion formulae for the operators $T_{p^\ell}$ with $p$ prime. 
It is the purpose of this paper to prove analogous results in the half-integral weight
setting. We also give an explicit formula for how operators $T_{p^{\ell}}$ commute
with the Shimura correspondence.

Let $k$, $N$ be positive integers with $k$ odd and $4 \mid N$.
Let $\chi$ be a Dirichlet character modulo $N$. 
We shall denote by $M_{k/2}(N,\chi)$ the space of
modular forms of weight $k/2$, level $N$ and character $\chi$,
and by $S_{k/2}(N,\chi)$ the subspace of cusp forms. 
We shall write $\TT_{k/2}$ for the Hecke algebra acting on these
spaces. For definitions we refer to Shimura's paper \cite{Shimura}.
It is well-known \cite[page 450]{Shimura} that $T_n=0$
for $n$ not square. 
\begin{thm}
Let $p$ be a prime and $\ell \ge 2$ be a positive integer. 
If $p \mid N$ then as Hecke operators in $\TT_{k/2}$,
\[
T_{p^{2 \ell}}=(T_{p^2})^\ell.
\]
If $p \nmid N$ then
\[
 T_{p^{2\ell+2}} = T_{p^2} T_{p^{2\ell}} - \chi(p^2)p^{k-2}T_{p^{2\ell-2}}
\] 
as Hecke operators in $\TT_{k/2}$.
\end{thm}

Let $N^\prime =N/2$. Let $t$ be a square-free positive integer. For $k \geq 3$,
Shimura defined the so-called \lq Shimura correspondence\rq
\[
\Sh_t: S_{k/2}(N,\chi) \rightarrow M_{k-1}(N^\prime,\chi^2).
\]
Here, $M_{k-1}(N^\prime,\chi^2)$ denotes the usual space of 
modular forms of integral weight $k-1$, level $N^\prime$ and character
$\chi^2$. It is well-known \cite{Kohnen} and 
\cite[page 53]{Ono}
that for $p \nmid t N$, 
\[
\Sh_t(T_{p^2} f)=T_p(\Sh_t(f))
\]
where $T_p$ is the usual integral weight Hecke operator. Our next theorem
shows that the same identity holds even when $p \mid tN$,
and also gives the precise relationship of how $T_{p^{2\ell}}$ commutes with $\Sh_t$ for all primes $p$,
and for $\ell \geq 2$.
\begin{thm} 
Let $p$ be a prime and let $f \in S_{k/2}(N,\chi)$. 
Let $t$ be a square-free positive integer. Then
\[
\Sh_t( T_{p^2} f)=T_p(\Sh_t(f)).
\]
Let $\ell \ge 2$. Then
\begin{enumerate}
\item[(a)] If $p \mid N$ then
\[
\Sh_t(T_{p^{2\ell}} f)=T_{p^\ell} (\Sh_t(f)).
\]
\item[(b)] If $p \nmid N$ then
\[
\Sh_t( T_{p^{2\ell}} f)= (T_{p^\ell} - \chi(p^2)p^{k-3}T_{p^{\ell-2}})( \Sh_t(f)).
\]
%where as before $T_{p^{2\ell}} \in \TT_{k/2}$ and $T_{p^\ell}$, $T_{p^{\ell-2}} \in \TT_{k-1}$.
\end{enumerate}
\end{thm}

It is well-known \cite[page 478]{Shimura} that the space $S_{3/2}(N,\chi)$ contains
single-variable theta-series and we denote by $S_0(N,\chi)$ the subspace generated 
by these single-variable theta-series. 
The interesting part of the space $S_{3/2}(N,\chi)$ is the orthogonal complement of 
$S_0(N,\chi)$ with respect to the Petersson inner product. We will use the following notation:
\[
S_{k/2}^{\perp}(N,\chi):= 
\begin{cases}
S_0(N,\chi)^\perp & k=3 \\
S_{k/2}(N,\chi) & \text{for $k \ge 5$}.
\end{cases}
\]
The Hecke algebra $\TT_{k/2}$ preserves $S_{k/2}^{\perp}(N,\chi)$.
We denote
\[
\TT_{k/2}^\perp=
\left\{
T\vert_{S_{k/2}^\perp (N,\chi)} : T \in \TT_{k/2}
\right\};
\]
this is the restriction of the Hecke algebra to $S_{k/2}^\perp (N,\chi)$.
\begin{thm} 
Let $k$, $N$ be positive integers with $k\geq 3$ odd,
and $4 \mid N$. Let $\chi$ be a Dirichlet character
modulo $N$. Let $N^\prime=N/2$. Write
\[
m= {N^\prime}^2 \prod_{p \mid N^\prime} 
\left( 1- \frac{1}{p^2} \right),
\qquad 
R= \frac{(k-1)m}{12}-\frac{m-1}{N^\prime}.
\]
Then $T_{i^2}$ for $i \leq R$
generate $\TT_{k/2}^\perp$ as a $\Z[\zeta_{\varphi(N)}]$-module. 
In particular the set of operators $T_{p^2}$ for primes $p \leq R$ forms a 
generating set as an algebra.
%Moreover, $f \in S_{k/2}(N,\chi)$
%is an eigenform for all Hecke operators if and only if it is
%an eigenform for $T_{p^2}$ for $p\leq R$.
%\end{cor}
%\begin{cor}\label{cor:sturm}
%With the same hypothesis as in the above corollary, further suppose that 
Moreover, if 
$\chi$ is a quadratic character, then the same result holds as above with 
\[
m = N^\prime \prod_{p \mid N^\prime} \left( 1+ \frac{1}{p} \right),
\qquad 
R= \frac{(k-1)m}{12}-\frac{m-1}{N^\prime}.
\]
\end{thm}

\section{Hecke operators}

\subsection{Integral weight Hecke operators}

\begin{comment}We recall the following definition of Hecke operators on the space of modular forms of integral weight 
in terms of double cosets.

Let $\alpha \in \GL_2^+(\R)$ such that $\Gamma_0(N)$ and $\alpha^{-1}\Gamma_0(N)\alpha$ are commensurable~
\footnote{Let $G$ be any group and $\Gamma$ and $\Gamma^\prime$ be two subgroups of $G$. We say that 
$\Gamma$ and $\Gamma^\prime$ are {\it commensurable} if 
\[[\Gamma:\Gamma\cap\Gamma^\prime]< \infty \qquad \text{and} \qquad [\Gamma^\prime: \Gamma\cap\Gamma^\prime]< \infty. \]}. 
Let $n$ be a positive integer. Then for any $f\in M_k(N)$ we have the following linear operators.
\begin{enumerate}
\item [(i)] \[f|[\Gamma_0(N)\alpha\Gamma_0(N)]_k(z) := \mathrm{det}(\alpha)^{k/2-1}\sum_{\nu=1}^{d}f|[\alpha_\nu]_k(z),
\] 
where $\alpha_\nu$ runs through the right cosets of $\Gamma_0(N)$ in $\Gamma_0(N)\alpha\Gamma_0(N)$.
\item[(ii)] \[T_n(f) := \sum f|[\Gamma_0(N)\alpha\Gamma_0(N)]_k, \]
where the sum is over all $\alpha =\left[ \begin{smallmatrix}l&0\\0&m\end{smallmatrix} \right]$ with $\ell,m$ 
positive integers, $\ell\mid m$, $(l,N)=1$ and $\ellm=n$.
\item[(iii)] If $(n,N)=1$, then
\[T_{(n,n)}(f) := f|[\Gamma_0(N)\left[ \begin{smallmatrix}n&0\\0&n\end{smallmatrix} \right]\Gamma_0(N)]_k.\]
\end{enumerate}
The operators $T_n$ and $T_{(n,n)}$ are called the {\it Hecke operators}.

The Hecke operators so defined preserve the cusp forms and one can similarly define the 
Hecke operators on the space of modular forms with characters. 
\end{comment}
Let $k$, $N$ be positive integers and $\chi$ be a Dirichlet character modulo $N$.
Let $M_k(N,\chi)$ be the space of modular forms of weight $k$, level $N$ and character $\chi$ and 
$S_k(N,\chi)$ be the subspace of cusp forms.
Recall that given a positive integer $n$ one can define Hecke operators $T_n$ and $T_{(n,n)}$ (when $(n,N)=1$) 
acting on the space $M_k(N,\chi)$ that also preserve $S_k(N,\chi)$. 

The following proposition lists the important properties of these Hecke operators.
\begin{prop}\label{prop:Heckeprop1} \begin{enumerate}
\item[(a)] If $(m,n)=1$, then $T_{mn} = T_mT_n$. 
\item[(b)] If $p$ is a prime dividing $N$, then $T_{p^e} = T_p^e$ for any positive integer $e$.
\item[(c)] If $p$ is a prime such that $(p,N)=1$, then for any positive integer $e$,
$ T_{p^{e+1}} = T_p T_{p^e} - p T_{(p,p)}T_{p^{e-1}}$ where for $f \in M_k(N,\chi)$ the action 
of $T_{(p,p)}$ can be explicitly expressed as $T_{(p,p)}(f) = p^{k-2}\chi(p) f$.
\end{enumerate}
\end{prop}
\begin{proof}
 See \cite[Lemma 4.5.7]{Miyake} and \cite[Pages 142-143]{Miyake}.
\end{proof}

The Hecke algebra on $M_k(N,\chi)$, 
which we denote by $\TT_k$ is an algebra over $\Z$ generated by $T_p$, $T_{(p,p)}$ and $T_q$ 
where $p$, $q$ varies over primes with $p \nmid N$ and $q \mid N$.
We can write the action of Hecke operators in terms of $q$-expansions.
\begin{prop}\label{prop:qheckeint}
Let $f$ be a modular form in $M_k(N,\chi)$ with $q$-expansion $f(z)=\sum_{n=0}^{\infty} a_n q^n$. 
Then $T_p(f)(z) = \sum_{n=0}^{\infty} b_n q^n$ where, 
\[b_n = a_{pn} +\chi(p)p^{k-1}a_{n/p}.\]
Here we take $a_{n/p} =0$ if $p \nmid n$.
\end{prop}
\begin{proof}
See \cite[Lemma 4.5.14]{Miyake}.
\end{proof}

\subsection{Half-integral weight forms and Hecke operators}\label{subsection:halfint}
Let $G$ be the group consisting of all ordered pairs $(\alpha,\phi(z))$, where 
$\alpha=\left[ \begin{smallmatrix}a&b\\c&d\end{smallmatrix} \right] \in \GL_2^+(\Q)$ 
and $\phi(z)$ is a holomorphic function on $\Hup$ satisfying
\[
 \phi(z)^2 = t \frac{cz+d}{\sqrt{\text{det}\alpha}}
\]for some $t \in \{\pm1\}$, with the group law defined by
\[
(\alpha,\phi(z))\cdot(\beta,\psi(z)) = (\alpha\beta,\phi(\beta z)\psi(z)). 
\]

Let $P:G\rightarrow\GL_2^+(\Q)$ be the homomorphism given by the projection map onto the first coordinate. 
Let $k$ be positive odd integer. The group $G$ acts on the space of complex valued functions on $\Hup$ by 
$f|[\xi]_{k/2}(z) := f(\alpha z)\phi(z)^{-k}$,
where $\xi=(\alpha,\phi(z)) \in G$ and $f:\Hup\rightarrow\C$.

Let $N$ be a positive integer with $4\mid N$. Then for $\gamma = \left[ \begin{smallmatrix}a&b\\c&d\end{smallmatrix} \right]
\in \Gamma_0(N)$ define
\[j(\gamma,z):=\kro{c}{d} \epsilon_d^{-1}\sqrt{cz+d},\qquad
\Delta_0(N):=\{\widetilde{\gamma}:=(\gamma,j(\gamma,z))|\gamma \in \Gamma_0(N)\}.\]
Note that $\Delta_0(N)$ is a subgroup of $G$. The map $L:\Gamma_0(4)\rightarrow G$ given by 
$\gamma\mapsto\widetilde{\gamma}$ defines an isomorphism onto $\Delta_0(4)$ and 
is mutually inverse to $P$. Denote by $\Delta_1(N)$ and $\Delta(N)$ respectively the images
of $\Gamma_1(N)$ and $\Gamma(N)$. 

Let $\chi$ be a Dirichlet character modulo $N$ and $M_{k/2}(N,\chi)$ and $S_{k/2}(N,\chi)$ be the 
spaces of modular forms and cusp forms of weight $k/2$, level $N$ and character $\chi$. 
The space $M_{k/2}(N,\chi)$ is $\{0\}$ unless $\chi$ is even, so henceforth we will be assuming 
$\chi$ to be even.
\begin{comment}
A holomorphic function $f$ on $\Hup$ is a {\it modular form of weight $k/2$} for $\Delta_1(N)$
if $f$ satisfies $f|[\widetilde{\gamma}]_{k/2} = f$ for all $\gamma \in \Gamma_1(N)$ and is 
holomorphic at all the cusps of $\Gamma_1(N)$, $f$ is called a cusp form if it vanishes
on all cusps. We denote such a space of modular forms 
by $M_{k/2}(\Gamma_1(N))$ and the subspace of cusp forms by $S_{k/2}(\Gamma_1(N))$. 
Further, let $M_{k/2}(N,\chi)$ 
(respectively $S_{k/2}(N,\chi)$) be the subspace of $M_{k/2}(\Gamma_1(N))$ 
(respectively $S_{k/2}(\Gamma_1(N))$) consisting of all elements $f$ such that
$f|[\widetilde{\gamma}]_{k/2} = \chi(d)f$ for all
$\gamma = \left[ \begin{smallmatrix}a&b\\c&d\end{smallmatrix} \right] \in \Gamma_0(N).$
\end{comment}
As in the integral weight case one can define the Hecke operators on the spaces
$M_{k/2}(N,\chi)$ and $S_{k/2}(N,\chi)$.

Let $\xi$ be an element of $G$ such that $\Delta_1(N)$ and $\xi^{-1}\Delta_1(N)\xi$ are commensurable.
Define an operator $|[\Delta_1(N)\xi\Delta_1(N)]_{k/2}$ on $M_{k/2}(\Gamma_1(N))$ by
\[f|[\Delta_1(N)\xi\Delta_1(N)]_{k/2} = \text{det}(\xi)^{k/4-1}\sum_{\nu}f|[\xi_\nu]_{k/2} \]
where $\Delta_1(N)\xi\Delta_1(N)=\bigcup_\nu\Delta_1(N)\xi_\nu$.

Now suppose $m$ is a positive integer and $\alpha = \left[ \begin{smallmatrix}1&0\\0&m\end{smallmatrix} \right]$, 
$\xi=(\alpha, m^{1/4})$. Then the Hecke operator $T_m$ is defined as the restriction of 
$|[\Delta_1(N)\xi\Delta_1(N)]_{k/2}$ to $M_{k/2}(N,\chi)$. It is to be noted that by 
\cite [Proposition 1.0] {Shimura}, if $m$ is not a square and $(m,N)=1$ then 
$|[\Delta_1(N)\xi\Delta_1(N)]_{k/2}$ is the zero operator. 
So we assume that $m = n^2$ for a positive integer $n$. Shimura writes the 
Hecke operator $T_{n^2}$ as 
\[T_{n^2}(f) := n^{\frac{k}{2}-2}\sum_\nu \chi(a_\nu)f|[\xi_\nu]_{k/2},\]
where $\xi_\nu$ are the right coset representatives of $\Delta_0(N)$ in 
$\Delta_0(N)\xi\Delta_0(N)$ such that 
$P(\xi_\nu)= \left[ \begin{smallmatrix}a_\nu&\ast\\\ast&\ast\end{smallmatrix} \right]$.
We have the following theorem.
\begin{thm}\label{thm:qheckehalf} (Shimura)
Let $f(z) = \sum_{n=0}^{\infty} a_n q^n \in M_{k/2}(N,\chi)$. Then $T_{p^2}(f)(z) = \sum_{n=0}^{\infty} b_n q^n$ where, 
\[b_n = a_{p^2n} + \chi(p){\kro{-1}{p}}^{\lambda}\kro{n}{p}p^{\lambda-1}a_n + \chi(p^2)p^{k-2}a_{n/p^2},\]
and $\lambda=(k-1)/2$ and $a_{n/p^2}=0$ whenever $p^2\nmid n$.
\end{thm}
\begin{proof}
See \cite[Theorem 1.7] {Shimura}.
\end{proof}

\section{Shimura Correspondence}\label{section:shimuralifts}

For this section fix positive integers $k$, $N$ with $k\geq 3$ odd and $4 \mid N$.
Let $\chi$ be an even Dirichlet character of modulus $N$.
Let $N^\prime=N/2$. We recall Shimura's Theorem.

\begin{thm} (Shimura) \label{thm:shimuracor}
Let $\lambda=(k-1)/2$. Let $f(z)=\sum_{n=1}^\infty a_n q^n \in S_{k/2}(N,\chi)$.
Let $t$ be a square-free integer and let $\psi_t$ be the
Dirichlet character modulo $tN$ defined by
\[
\psi_t(m)=\chi(m) \kro{-1}{m}^\lambda  \kro{t}{m}.
\]
Let $A_t(n)$ be the complex numbers defined by
\begin{equation}\label{eqn:shiexp}
\sum_{n=1}^\infty A_t(n) n^{-s} = \left( \sum_{i=1}^\infty \psi_t(i) i^{\lambda-1-s} \right)
\left( \sum_{j=1}^\infty a_{t j^2} j^{-s} \right).
\end{equation}
Let $\Sh_t(f)(z)=\sum_{n=1}^\infty A_t(n) q^n$.
Then
\begin{enumerate}
\item[(i)] $\Sh_t(f) \in M_{k-1}(N^\prime,\chi^2)$.
\item[(ii)] If $k \geq 5$ then $\Sh_t(f)$ is a cusp form.
\item[(iii)] If $k =3 $ and $f \in S_{3/2}^\perp(N,\chi)$
then $\Sh_t(f)$ is a cusp form. 
%\item[(iv)] Suppose $f$ is an eigenform for $T_{p^2}$ for all primes $p$ and let $T_{p^2}f = \lambda_p f$. 
%Then $\sum_{n=1}^\infty A_0(n) q^n \in M_{k-1}(N^\prime,\chi^2)$  where $A_0(n)$ is defined by 
%\begin{equation}\label{eqn:shiexp1}
%\sum_{n=1}^\infty A_0(n) n^{-s} = \prod_p (1 - \lambda_p p^{-s} + \chi(p)^2 p^{k-2-2s})^{-1}.
%\end{equation}
%In fact if $a_t \ne 0$ then $\Sh_t(f)/a_t = \sum_{n=1}^\infty A_0(n) q^n$.
\end{enumerate} 
\end{thm}
\begin{proof}
For (i) and (ii) see \cite[Section 3, Main Theorem]{Shimura}, for the rest see
\cite[Theorem 3.14]{Ono}. In particular, the fact that $N^\prime=N/2$ was proved by Niwa \cite[Section 3]{Niwa}.
\end{proof}
The form $\Sh_t(f)$ is called the {\em Shimura lift of $f$
corresponding to $t$}.
The following is clear from Equation\eqref{eqn:shiexp}.
\begin{lem}\label{lem:shlin}
The Shimura lift $\Sh_t$ is linear.
\end{lem}

\begin{lem}\label{lem:shker}
If $\Sh_t(f)=0$ for all positive square-free integers $t$ then
$f=0$.
\end{lem}
\begin{proof}
By Equation~\eqref{eqn:shiexp} we know that $a_{tj^2}=0$ for
all positive square-free integers $t$ and all positive integers $j$.
Then $a_n=0$ for all $n$.
\end{proof}

In Ono's book \cite[Chapter 3, Corollary 3.16]{Ono} and several other places 
\cite{Kohnen} we find the following result stated without proof.
\begin{prop}
Suppose $f \in S_{k/2}(N,\chi)$.  Let $t$ be a square-free positive integer. If $p \nmid tN$
is a prime then
\[
\Sh_t( T_{p^2} f)=T_{p} \Sh_t(f). 
\]
\end{prop}
Here $T_{p^2}$ is the Hecke operator in $\TT_{k/2}$ and $T_p$ is the 
Hecke operator in $\TT_{k-1}$. For what follows we shall need the following strengthening 
of this result.
\begin{prop}\label{prop:commute}
Suppose $f \in S_{k/2}(N,\chi)$ and $t$ a square-free
positive integer.  If $p$
is a prime then
\[
\Sh_t( T_{p^2} f)=T_{p} \Sh_t(f). 
\]
\end{prop}
We do not know why the above references impose 
the condition $p \nmid tN$. We shall give a careful proof
that does not use this assumption.
\begin{proof}[Proof of Proposition~\ref{prop:commute}]
The proof uses the explicit formulae for Hecke operators in terms of $q$-expansions.
As in Shimura's Theorem above, write $f(z)=\sum_{n=1}^\infty a_n q^n$. Fix $t$ to be a 
positive square-free integer.
To simplify notation, we shall write $A_n$ for $A_t(n)$. Thus we have the relation
\[
\sum_{n=1}^\infty A_n n^{-s} = \left( \sum_{i=1}^\infty \psi_t(i) i^{\lambda-1-s} \right)
\left( \sum_{j=1}^\infty a_{t j^2} j^{-s} \right).
\] 
We may rewrite this as
\begin{equation}\label{eqn:An}
A_n= \sum_{ij=n} \psi_t(i) i^{\lambda-1} a_{tj^2}.
\end{equation} 
Let
\[
T_{p^2}(f)(z)=\sum_{n=1}^\infty b_n q^n.
\]
Then using Theorem~\ref{thm:qheckehalf} we get,
\begin{equation}\label{eqn:bn}
b_n=a_{p^2 n} +\psi_1(p) \kro{n}{p} p^{\lambda-1} a_n + 
\chi^2(p) p^{k-2} a_{n/p^2}.
\end{equation}
The reader will recall that if $n/p^2$ is not an integer then we take $a_{n/p^2}=0$.
%We note that $\psi_1^2=\chi^2$; we shall make use of this fact several times in proof.

Let $g=\Sh_t(f)(z)=\sum_{n=1}^\infty A_n q^n$. Write
\[
T_p(g)(z)=\sum_{n=1}^\infty B_n q^n.
\]
Let 
\[
\Sh_t(T_{p^2} f)(z)=\sum_{n=1}^\infty C_n q^n.
\]
To prove the proposition, it is enough to show that $B_n=C_n$ for all $n$. We shall do this
by direct calculation, expressing both $B_n$ and $C_n$ in terms of the $a_i$.

Since $g(z)=\sum A_n q^n \in M_{k-1}(N^\prime,\chi^2)$ and $T_p(g)(z)=\sum B_n q^n$ we know 
by Proposition \ref{prop:qheckeint} that
\[
B_n= A_{pn}+ \chi^2(p) p^{k-2} A_{n/p}.
\]
Substituting from \eqref{eqn:An} we have
\begin{equation}\label{eqn:Bn}
B_n=\sum_{ij=pn} \psi_t(i) i^{\lambda-1} a_{tj^2} +
\sum_{ij=n/p} \chi^2(p) \psi_t(i) p^{k-2} i^{\lambda-1} a_{tj^2};
\end{equation}
here the second sum is understood to vanish if $p \nmid n$.

Recall $T_{p^2}f(z)=\sum b_n q^n$ and $\Sh_t(T_{p^2} f)(z)=\sum C_n q^n$.
Hence by \eqref{eqn:An} we have
\[
C_n=\sum_{ij=n} \psi_t(i) i^{\lambda-1} b_{tj^2}.
\]
Using \eqref{eqn:bn} we obtain
\[
C_n= \sum_{ij=n} \psi_t(i) i^{\lambda-1} 
\left(
a_{p^2 tj^2} +\psi_1(p) \kro{tj^2}{p} p^{\lambda-1} a_{tj^2}+
\chi^2(p) p^{k-2} a_{tj^2/p^2}
\right).
\]
Note that $\psi_1(p)\kro{tj^2}{p} = \psi_t(p)\kro{j^2}{p}$.
So we can rewrite $C_n$ as
\begin{equation}\label{eqn:Cn}
C_n= \sum_{ij=n} \psi_t(i) i^{\lambda-1} 
\left(
a_{p^2 tj^2} +\psi_t(p) \kro{j^2}{p} p^{\lambda-1} a_{tj^2}+
\chi^2(p) p^{k-2} a_{tj^2/p^2}
\right).
\end{equation}
Note that the Legendre symbol here is $1$ unless of course $p \mid j$
in which case it is $0$.
Moreover $a_{tj^2/p^2} = 0$ whenever $p \nmid j$;
this is because $t$ is square-free. 

We consider the following two cases.

\noindent {\bf Case $p \nmid n$}.
In this case the formulae for $B_n$ and $C_n$ simplify as follows.
\begin{equation*}
\begin{split}
B_n &=\sum_{ij=pn} \psi_t(i) i^{\lambda-1} a_{tj^2}\\
&= \sum_{ij=n} \psi_t(pi) (pi)^{\lambda-1} a_{tj^2}+ \psi_t(i) i^{\lambda-1} a_{tp^2 j^2}\\
&= \sum_{ij=n} \psi_t(i) i^{\lambda-1} (a_{tp^2 j^2} +\psi_t(p) p^{\lambda -1} a_{tj^2})\\
&= C_n.
\end{split}
\end{equation*}

\noindent {\bf Case $p \mid n$}. Write $n=p^r m$ where $r \geq 1$ and $p \nmid m$. We rewrite
\eqref{eqn:Bn} as follows.
\begin{multline*}
B_n = 
\sum_{j \mid p^{r+1} m}  \psi_t(p^{r+1} m/j) (p^{r+1} m/j)^{\lambda-1} a_{tj^2} \\
+\sum_{j \mid p^{r-1} m} \chi^2(p) \psi_t(p^{r-1} m/j) p^{k-2} (p^{r-1} m/j)^{\lambda-1} a_{tj^2}.
\end{multline*}
This maybe re-expressed as $B_n=B_n^{(1)}+B_n^{(2)}$ where
\[
B_n^{(1)} = 
\sum_{u=0}^{r+1}\sum_{k \mid m}  \psi_t(p^{r+1-u} m/k) (p^{r+1-u} m/k)^{\lambda-1} a_{tp^{2u} k^2}
\]
and
\[
B_n^{(2)}= 
\sum_{u=0}^{r-1}\sum_{k \mid m} \chi^2(p) \psi_t(p^{r-1-u} m/k) p^{k-2} (p^{r-1-u} m/k)^{\lambda-1} a_{tp^{2u} k^2}.
\]

Moreover, we can rewrite \eqref{eqn:Cn} as follows.
\[
C_n= \sum_{j \mid p^r m} \psi_t(p^r m/j) (p^r m/j)^{\lambda-1}
\left(
a_{p^2 t j^2}+\psi_t(p) \kro{j^2}{p} p^{\lambda-1} a_{tj^2} 
+\chi^2(p) p^{k-2} a_{tj^2/p^2} \right).
\]
Thus we can write $C_n=C_n^{(1)}+C_n^{(2)}+C_n^{(3)}$ where
\[
C_n^{(1)}=
\sum_{u=0}^r \sum_{k \mid m}
\psi_t(p^{r-u} m/k) (p^{r-u} m/k)^{\lambda-1} a_{t p^{2u+2} k^2},
\]
and
\[
C_n^{(2)}=
\sum_{k \mid m}
\psi_t(p^{r+1} m/k) (p^{r+1} m/k)^{\lambda-1}  a_{t k^2},
\]
and
\[
C_n^{(3)}=
\sum_{u=1}^r \sum_{k \mid m}
\chi^2(p) \psi_t(p^{r-u} m/k) (p^{r-u} m/k)^{\lambda-1}  p^{k-2} a_{t p^{2u-2} k^2}.
\]
It is clear that 
$B_n^{(2)}=C_n^{(3)}$, and also that
$B_n^{(1)}=C_n^{(1)}+C_n^{(2)}$; here $C_n^{(2)}$ corresponds
to the $u=0$ terms in $B_n^{(1)}$.
Thus $B_n=C_n$ completing the proof.
\end{proof}

\section{Recursion Formula for the Hecke Operators $T_{p^{2\ell}}$}

We keep the notation as in the previous section.   
Let $\ell$ be a positive integer and $p$ be a prime. In this section we are 
interested in the action of the Hecke operator $T_{p^{2\ell}}$ on the space $M_{k/2}(N,\chi)$. 
In the case $p \mid N$ we have the following easy lemma. 
\begin{lem}
Let $\ell$ be a positive integer and $p$ be a prime dividing $N$. 
Let $t$ be a square-free positive integer. Then
\begin{enumerate}
\item[(i)] $T_{p^{2\ell}} = (T_{p^2})^\ell$.
\item[(ii)] $\Sh_t(T_{p^{2\ell}} f) = T_{p^\ell}(\Sh_t (f))$ for $f \in S_{k/2}(N,\chi)$.
\end{enumerate}
In the above statements $T_{p^{2\ell}} \in \TT_{k/2}$ and $T_{p^\ell} \in \TT_{k-1}$.
\end{lem}
\begin{proof}
Let $f=\sum_{n=0}^{\infty}a_n q^n \in M_{k/2}(N,\chi)$. It follows using ~\cite[Proposition 1.5]{Shimura} that
$T_{p^{2\ell}}(f) = \sum_{n=1}^{\infty}a_{np^{2\ell}} q^n$. Now part $(i)$ follows using 
Theorem~\ref{thm:qheckehalf}. Part $(ii)$ follows by using Proposition~\ref{prop:commute} and 
part $(b)$ of Proposition~\ref{prop:Heckeprop1} since $p \mid N^\prime$.
\end{proof}

We will assume that $p \nmid N$ for the rest of this section.
The main aim of this section is to prove the following result. 

\begin{thm}\label{thm:recur}
Let $p \nmid N$ be a prime and $\ell \ge 2$ be a positive integer. Then 
\[
 T_{p^{2\ell+2}} = T_{p^2} T_{p^{2\ell}} - \chi(p^2)p^{k-2}T_{p^{2\ell-2}}
\] 
as Hecke operators in $\TT_{k/2}$.
\end{thm}
It is to be noted that for $l=1$ the above relation does not hold. One can check directly that in $\TT_{k/2}$,
\[T_{p^4} = (T_{p^2})^2 - \chi(p^2)(p^{k-3} + p^{k-2}).\]

We need the following lemma on Gauss sums which can be easily deduced from \cite[Lemma 3.1.3]{Miyake}:
\begin{lem}\label{lem:Gauss}
Let $p$ be a prime and $n$, $\alpha$ be a given positive integer. Then
\begin{enumerate}
\item[(i)]$\sum_{m=0}^{p^{\alpha}-1} \kro{m}{p} e^{\frac{2\pi imn}{p^\alpha}} =
\begin{cases}
0  &   \text{if $p^{\alpha-1} \nmid n$}\\
p^{\alpha-1}\kro{n^\prime}{p}\epsilon_p \sqrt{p} &   \text{if $n=p^{\alpha-1}n^\prime$}.
\end{cases}$
\item[(ii)]
$\sum_{m=0}^{p^{\alpha}-1} e^{\frac{2\pi imn}{p^\alpha}} =
\begin{cases}
0  &   p^{\alpha} \nmid n\\
p^{\alpha} &   p^{\alpha} \mid n.
\end{cases}$
\end{enumerate}
\end{lem}

\begin{proof}[Proof of Theorem~\ref{thm:recur}]
Let $f \in M_{k/2}(N,\chi)$. Let $\alpha = [\begin{smallmatrix} 1 & 0\\ 0& p^{2\ell} \end{smallmatrix}]$, 
$\xi = (\alpha, p^{\ell/2})$. Using \cite[Lemma 4.5.6]{Miyake} we know that 
\[
 \Gamma_0(N) \alpha \Gamma_0(N) = \bigcup_{\nu,m} \Gamma_0 \alpha_{\nu,m}, \qquad 
\alpha_{\nu,m} = \left[\begin{matrix} p^{2\ell-\nu} & m\\ 0& p^{\nu} \end{matrix}\right]
\]
where $0 \le \nu \le 2\ell$, $0 \le m < p^{\nu}$ and $\gcd(m, p^{\nu}, p^{2\ell-\nu})=1$. 
Let $G$ be the group defined as above.
Let $\xi_{\nu,m} \in G$ be given by 
\[
\xi_{\nu,m} = 
\begin{cases} 
\left(\alpha_{\nu,m},\ p^{\frac{-2\ell + 2\nu}{4}} \epsilon_p^{-1}\kro{-m}{p} \right) &  \text{if $\nu$ is odd}\\
(\alpha_{\nu,m},\ p^{\frac{-2\ell + 2\nu}{4}}) &  \text{if $\nu$ is even}.
\end{cases}
\]
One can verify that $\xi_{\nu,m}$ with $\nu$ and $m$ varying as above form a set of right coset 
representatives of $\Delta_0(N)$ in $\Delta_0(N)\xi\Delta_0(N)$ (see \cite[Proposition 1.1]{Shimura}).
Then we know by definition of $T_{p^{2\ell}}$ (see Subsection~\ref{subsection:halfint}) that 
\begin{equation}\label{eqn:A}
T_{p^{2\ell}} f = (p^{2\ell})^{\frac{k}{4}-1} \left(
A_0+ A_{2\ell} 
+
\sum_{\nu=1}^{2\ell-1} A_{\nu}
\right), 
\end{equation}
where 
\[
 A_{\nu} = \sum_{\substack{m=0 \\ (m,p)=1}}^{p^\nu -1}
\chi(p^{2\ell-\nu}) f|[\xi_{\nu,m}]_{k/2}, \ A_{2\ell} = \sum_{m=0}^{p^{2\ell} -1} f|[\xi_{2\ell,m}]_{k/2}, \ 
A_0 = \chi(p^{2\ell}) f|[\xi_{0,0}]_{k/2}.
\]
\begin{comment}
\begin{equation*}
\begin{split}
T_{p^{2\ell}} f &= (p^{2\ell})^{\frac{k}{4}-1} 
\left( \sum_{\nu=1}^{2\ell-1} \sum_{\substack{m=0 \\ \gcd(m,p)=1}}^{p^\nu -1}
\chi(p^{2\ell-\nu}) f|[\xi_{\nu,m}]_{k/2} + 
%\sum_{\substack{\nu=2\\ \nu \text{odd}}}^{2\ell-2} \sum_{\substack{m=0 \\ \gcd(m,p)=1}}^{p^\nu -1}
%\chi(p^{2\ell-\nu}) f|[\xi_{\nu,m}]_{k/2} + 
\sum_{m=0}^{p^{2\ell} -1} f|[\xi_{2\ell,m}]_{k/2} +
\chi(p^{2\ell}) f|[\xi_{0,0}]_{k/2}\right)\\
&= (p^{2\ell})^{\frac{k}{4}-1} (\sum_{\nu=1}^{2\ell-1} A_{\nu} + A_{2\ell} + A_0)
\end{split}
\end{equation*}
\end{comment}
Applying $T_{p^2}$ to Equation~\eqref{eqn:A} we obtain 
\begin{equation}\label{eqn:B}
\begin{split}
T_{p^2} T_{p^{2\ell}} f &= (p^{2\ell})^{\frac{k}{4}-1} \left(\sum_{\nu=1}^{2\ell-1} T_{p^2} A_{\nu} + 
T_{p^2} A_{2\ell} + T_{p^2}A_0\right)\\
&= (p^{2\ell+2})^{\frac{k}{4}-1} \left(\sum_{\nu=1}^{2\ell-1} B_{\nu} + B_{2\ell} + B_0\right),
\end{split}
\end{equation}
where for $\nu$ with $0 \le \nu \le 2\ell-2$ we have 
\begin{equation*}
\begin{split}
B_\nu &= \chi(p^{2\ell-\nu +2}) \sum_{\substack{m=0 \\ (m,p)=1}}^{p^\nu -1} 
f|[([\begin{smallmatrix}p^{2\ell-\nu +2} & m \\ 0 & p^{\nu} \end{smallmatrix}],\  
p^{\frac{-2\ell + 2\nu -2}{4}} r_{\nu,m})]_{k/2} \\
& + \chi(p^{2\ell-\nu +1}) \sum_{m^\prime=1}^{p-1}\sum_{\substack{m=0 \\ (m,p)=1}}^{p^\nu -1} 
f|[([\begin{smallmatrix}p^{2\ell-\nu +1} & p^{2\ell-\nu}m^\prime + mp \\ 0 & p^{\nu+1} \end{smallmatrix}],\  
p^{\frac{-2\ell + 2\nu}{4}} s_{\nu,m,m^\prime})]_{k/2} \\
& + \chi(p^{2\ell-\nu}) \sum_{m^\prime=0}^{p^2-1}\sum_{\substack{m=0 \\ (m,p)=1}}^{p^\nu -1} 
f|[([\begin{smallmatrix}p^{2\ell-\nu} & p^{2\ell-\nu}m^\prime + mp^2 \\ 0 & p^{\nu+2} \end{smallmatrix}],\  
p^{\frac{-2\ell + 2\nu +2}{4}} r_{\nu,m})]_{k/2},
\end{split}
\end{equation*}
where
\[
r_{\nu,m} =\begin{cases}
\epsilon_p^{-1}\kro{-m}{p} & \text{$\nu$ odd}\\
1 & \text{$\nu$ even}
\end{cases}, \qquad
s_{\nu,m,m^\prime} =\begin{cases}
\epsilon_p^{-2}\kro{mm^\prime}{p} & \text{$\nu$ odd}\\
\epsilon_p^{-1}\kro{-m^\prime}{p} & \text{$\nu$ even},
\end{cases}
\]
and $B_{2\ell}$ has the same expression as above with $\nu=2\ell$ but without any coprimality condition 
on $m$, that is, we do not have $(m,p)=1$ in the above terms while writing the expression for $B_{2\ell}$.
\begin{comment}
for $\nu$ even with $0 \le \nu \le 2\ell-2$ we have 
\begin{equation*}
\begin{split}
B_\nu = & \chi(p^{2\ell-\nu +2}) \sum_{\substack{m=0 \\ (m,p)=1}}^{p^\nu -1} 
f|[([\begin{smallmatrix}p^{2\ell-\nu +2} & m \\ 0 & p^{\nu} \end{smallmatrix}], 
p^{\frac{-2\ell + 2\nu -2}{4}})]_{k/2} + \\
& \chi(p^{2\ell-\nu +1}) \sum_{m^\prime=1}^{p-1}\sum_{\substack{m=0 \\ (m,p)=1}}^{p^\nu -1} 
f|[([\begin{smallmatrix}p^{2\ell-\nu +1} & p^{2\ell-\nu}m^\prime + mp \\ 0 & p^{\nu+1} \end{smallmatrix}], 
p^{\frac{-2\ell + 2\nu}{4}} \epsilon_p^{-1}\kro{-m^\prime}{p})]_{k/2} + \\
& \chi(p^{2\ell-\nu}) \sum_{m^\prime=0}^{p^2-1}\sum_{\substack{m=0 \\ (m,p)=1}}^{p^\nu -1} 
f|[([\begin{smallmatrix}p^{2\ell-\nu} & p^{2\ell-\nu}m^\prime + mp^2 \\ 0 & p^{\nu+2} \end{smallmatrix}], 
p^{\frac{-2\ell + 2\nu +2}{4}})]_{k/2},
\end{split}
\end{equation*}
\end{comment}

We express $T_{p^{2\ell+2}}f$ as in Equation~\eqref{eqn:A} and 
compare it with Equation~\eqref{eqn:B}. Ruling out some of the terms using Euclidean 
algorithm and rewriting the action of matrices (we will give an example of the working later) 
we obtain
\begin{equation}\label{eqn:AB}
(T_{p^{2\ell+2}}- T_{p^2} T_{p^{2\ell}})(f) = -(p^{2\ell+2})^{\frac{k}{4}-1} \left( S_0 + S_{2\ell} + 
\sum_{\nu=1}^{2\ell-1} ( D_\nu + E_\nu)\right)
\end{equation} where
\begin{equation*}
\begin{split}
S_0 &= \sum_{m^\prime=0}^{p^2 -1} \chi(p^{2\ell})
f|[([\begin{smallmatrix}p^{2\ell} & p^{2\ell}m^\prime \\ 0 & p^2 \end{smallmatrix}],\  
p^{\frac{-\ell+1}{2}})]_{k/2}\\
S_{2\ell} &= \sum_{\substack{m=0 \\ (m,p)\ne1}}^{p^{2\ell} -1} \chi(p^2)
f|[([\begin{smallmatrix}p^2 & m \\ 0 & p^{2\ell} \end{smallmatrix}],\ 
p^{\frac{\ell-1}{2}})]_{k/2} \\
D_\nu &= \chi(p^{2\ell-\nu})
\sum_{m^\prime=0}^{p^2-1}\sum_{\substack{m=0 \\ (m,p)=1}}^{p^\nu -1} 
f|[([\begin{smallmatrix}p^{2\ell-\nu} & p^{2\ell-\nu}m^\prime + mp^2 \\ 0 & p^{\nu+2} \end{smallmatrix}],\  
p^{\frac{-2\ell + 2\nu + 2}{4}} r_{\nu,m})]_{k/2}\\
E_\nu &= \chi(p^{2\ell-\nu +1})
\sum_{m^\prime=1}^{p-1}\sum_{\substack{m=0 \\ (m,p)=1}}^{p^\nu -1} 
f|[([\begin{smallmatrix}p^{2\ell-\nu +1} & p^{2\ell-\nu}m^\prime + mp \\ 0 & p^{\nu+1} \end{smallmatrix}],\  
p^{\frac{-2\ell + 2\nu}{4}} s_{\nu,m,m^\prime})]_{k/2}.
\end{split}
\end{equation*}
Further
\begin{equation}\label{eqn:C}
\chi(p^2)p^{k-2}T_{p^{2\ell-2}} f = p^2(p^{2\ell+2})^{\frac{k}{4}-1} \left(\sum_{\nu=1}^{2\ell-3} C_{\nu} + C_{2\ell-2} + C_0\right), 
\end{equation}
where for $\nu$ with $0\le \nu \le 2\ell-3$ we have
\[
C_{\nu} = \sum_{\substack{m=0 \\ (m,p)=1}}^{p^\nu -1}
\chi(p^{2\ell-\nu}) f|[([\begin{smallmatrix}p^{2\ell-\nu -2} & m \\ 0 & p^{\nu} \end{smallmatrix}],\  
p^{\frac{-2\ell + 2\nu + 2}{4}} r_{\nu,m})]_{k/2}\]
and $C_{2\ell-2}$ has the same expression as above with $\nu=2\ell-2$ but without the condition $(m,p)=1$ 
in the above sum. We first claim that the following relations hold:
\begin{enumerate}
\item[(i)] $ D_{\nu} = p^2 C_{\nu}$ for $1 \le \nu \le 2\ell-3, \quad \text{and} \quad S_0 = p^2C_0$.
\item[(ii)] $E_\nu = 0$ for $1 \le \nu \le 2\ell-2$.
%$ \sum_{\substack{\nu=1 \\ \nu \text{odd}}}^{2\ell-3} E_\nu = 0 \quad \text{and} \quad
%\sum_{\substack{\nu=2 \\ \nu \text{even}}}^{2\ell-2} E_\nu = 0.$
\end{enumerate}
We will only show the computation for part $(ii)$ for case $\nu$ odd. The rest of the claim follows by 
similar method. Fix an odd $\nu$ with $1 \le \nu \le 2\ell-3$. Fix $1 \le m^\prime \le p-1$. Then for 
each $m$ with $0 \le m \le p^{\nu}-1$ there exist unique $a$ and $b$ with $0 \le b \le p^{\nu}-1$ 
such that $m + p^{2\ell-\nu-1}m^\prime = ap^{\nu} + b$. Moreover $m \equiv b \pmod{p}$. Hence 
\[ (m,p)=1 \iff (b,p)=1, \qquad \kro{-m}{p} = \kro{-b}{p}.\]
We can rewrite $E_\nu$ as
\begin{equation*}
\begin{split}
&E_{\nu} = \chi(p^{2\ell-\nu +1})
\sum_{m^\prime=1}^{p-1}\sum_{\substack{m=0 \\ (m,p)=1}}^{p^\nu -1} 
f\left(\frac{p^{2\ell-\nu +1}z + p^{2\ell-\nu}m^\prime + mp} {p^{\nu+1}}\right)
\left(p^{\frac{-2\ell + 2\nu}{4}}\epsilon_p^{-2}\kro{mm^\prime}{p}\right)^{-k}\\
&= \chi(p^{2\ell-\nu +1})\epsilon_p^k\sum_{m^\prime=1}^{p-1}\kro{-m^\prime}{p}
\sum_{\substack{m=0 \\ (m,p)=1}}^{p^\nu -1}
f\left| \left[\left([\begin{smallmatrix}p^{2\ell-\nu} & p^{2\ell-\nu-1}m^\prime + m \\ 0 & p^{\nu} \end{smallmatrix}],\  
p^{\frac{-2\ell + 2\nu}{4}} \epsilon_p^{-1}\kro{-m}{p} \right)\right]_{k/2} \right.\\
&= \chi(p^{2\ell-\nu +1})\epsilon_p^k\sum_{m^\prime=1}^{p-1}\kro{-m^\prime}{p}
\sum_{\substack{b=0 \\ (b,p)=1}}^{p^\nu -1}
f \left|\left[\left([\begin{smallmatrix}p^{2\ell-\nu} & b \\ 0 & p^{\nu} \end{smallmatrix}],\  
p^{\frac{-2\ell + 2\nu}{4}} \epsilon_p^{-1}\kro{-b}{p}\right)\right]_{k/2} \right.\\
&= 0.
\end{split}
\end{equation*}
The second last equality follows since as elements of $G$ we have
\[
\left([\begin{smallmatrix}p^{2\ell-\nu} & p^{2\ell-\nu-1}m^\prime + m \\ 0 & p^{\nu} \end{smallmatrix}],\  
p^{\frac{-2\ell + 2\nu}{4}} \epsilon_p^{-1}\kro{-m}{p} \right)
= ([\begin{smallmatrix}1 & a \\ 0 & 1 \end{smallmatrix}],1) 
\cdot \left([\begin{smallmatrix}p^{2\ell-\nu} & b \\ 0 & p^{\nu} \end{smallmatrix}],\  
p^{\frac{-2\ell + 2\nu}{4}} \epsilon_p^{-1}\kro{-b}{p} \right).
\]
By working out similarly as above one can further see that 
\[p^2C_{2\ell-2} - D_{2\ell-2} = \chi(p^2) \sum_{m^\prime=0}^{p^2-1}
\sum_{\substack{m=0 \\ (m,p) \ne 1}}^{p^{2\ell-2} -1}
f|[([\begin{smallmatrix}p^2 & p^2m^\prime + mp^2 \\ 0 & p^{2\ell} \end{smallmatrix}],\  
p^{\frac{\ell-1}{2}})]_{k/2} =: F_{2\ell-2}.\]
Thus to prove the theorem we are left to show that 
\[F_{2\ell-2}-S_{2\ell}-E_{2\ell-1}-D_{2\ell-1} = 0.\]
We claim that $D_{2\ell-1} =0$ and $F_{2\ell-2}-S_{2\ell}-E_{2\ell-1}=0$ which proves the theorem.

We first show that $D_{2\ell-1} =0$. Let $f(z)= \sum_{n=0}^\infty a_n e(nz)$ where 
$e(nz)=e^{2\pi i nz}$. Rewriting $D_{2\ell-1}$ in terms of coefficients $a_n$ we obtain

\begin{equation*}
\begin{split}
&D_{2\ell-1} = \chi(p) p^{\frac{-\ell k}{2}}\epsilon_p^k\kro{-1}{p}
\sum_{m^\prime=0}^{p^2-1}\sum_{\substack{m=0 \\ (m,p)=1}}^{p^{2\ell-1} -1} 
\sum_{n=0}^\infty a_n e\left(\frac{npz + npm^\prime +nmp^2}{p^{2\ell+1}}\right)\kro{m}{p}\\
&= \chi(p) p^{\frac{-\ell k}{2}}\epsilon_p^k\kro{-1}{p} 
\sum_{n=0}^\infty a_n e\left(\frac{nz}{p^{2\ell}}\right)
\sum_{m^\prime=0}^{p^2-1}e\left(\frac{nm^\prime}{p^{2\ell}}\right)
\sum_{m=0}^{p^{2\ell-1}-1}e\left(\frac{nm}{p^{2\ell-1}}\right)\kro{m}{p}\\
&=  \chi(p) p^{\frac{-\ell k+4\ell-3}{2}}\epsilon_p^{k+1}\kro{-1}{p}
\sum_{\substack{n=0\\p^{2\ell-2}\mid n}}^\infty a_n e\left(\frac{nz}{p^{2\ell}}\right)
\kro{n/p^{2\ell-2}}{p}
\sum_{m^\prime=0}^{p^2-1}e\left(\frac{nm^\prime/p^{2\ell-2}}{p^2}\right)\\
&=0,
\end{split}
\end{equation*}
where last two equalities follows using Lemma~\ref{lem:Gauss} on Gauss sums. 
In order to prove the final claim we again use the coefficients method as above 
to obtain
\begin{equation*}
\begin{split}
F_{2\ell-2}-S_{2\ell} &= \chi(p^2)p^{\frac{(-\ell+1)k+4\ell-2}{2}}
\sum_{\substack{n=0\\p^{2\ell-2}\| n}}^\infty a_n e\left(\frac{nz}{p^{2\ell-2}}\right), \\
E_{2\ell-1} &= \chi(p^2)p^{\frac{(-\ell+1)k+4\ell-2}{2}} \epsilon_p^{2k+2}
\sum_{\substack{n=0\\p^{2\ell-2}\| n}}^\infty a_n e\left(\frac{nz}{p^{2\ell-2}}\right).
\end{split}
\end{equation*}
Now $\epsilon_p^{2k+2}=1$ since $2k+2 \equiv 0 \pmod{4}$. Hence we are done. 
\end{proof}

\begin{cor}
Let $p \nmid N$ be a prime and  $\ell\ge 2$. Let $f \in S_{k/2}(N,\chi)$. Then
\[
\Sh_t( T_{p^{2\ell}} f)= (T_{p^\ell} - \chi(p^2)p^{k-3}T_{p^{\ell-2}})( \Sh_t f),
\]
where as before $T_{p^{2\ell}} \in \TT_{k/2}$ and $T_{p^\ell}$, $T_{p^{\ell-2}} \in \TT_{k-1}$.
\end{cor}
\begin{proof}
We use induction on $\ell$. Recall from part $(c)$ of Proposition~\ref{prop:Heckeprop1} that for prime 
$p \nmid N$, we have
\begin{equation}\label{eqn:Heckeprop}
T_{p^{e+1}} (\Sh_t f) = (T_p T_{p^e}  - \chi(p^2)p^{k-2} T_{p^{e-1}})(\Sh_t f). 
\end{equation}
As we remarked earlier, for $l=2$ we have the following relation in $\TT_{k/2}$:
\[T_{p^4} = (T_{p^2})^2 - \chi(p^2)(p^{k-3} + p^{k-2}).\]
Hence we get
\begin{equation*}
\begin{split}
\Sh_t(T_{p^4} f) &= \Sh_t( (T_{p^2})^2 f) - \chi(p^2)(p^{k-3} + p^{k-2})(\Sh_t f)\\
&= ((T_p)^2 - \chi(p^2)p^{k-2})(\Sh_t f) - \chi(p^2)p^{k-3} (\Sh_t f) \\
&= ( T_{p^2} - \chi(p^2)p^{k-3} )(\Sh_t f).
\end{split}
\end{equation*}
Assume the statement holds for all $\ell \le e$. Then 
\begin{equation*}
\begin{split}
&\Sh_t( T_{p^{2e+2}} f) = \Sh_t ( T_{p^2} T_{p^{2e}} f)- \chi(p^2)p^{k-2} \Sh_t( T_{p^{2e-2}} f) \\
&= T_p (\Sh_t (T_{p^{2e}} f)) - \chi(p^2)p^{k-2} \Sh_t( T_{p^{2e-2}} f)\\
%&= T_p (T_{p^e} - \chi(p^2)p^{k-3}T_{p^{e-2}})( \Sh_t(f)) - 
%\chi(p^2)p^{k-2} (T_{p^{e-1}} - \chi(p^2)p^{k-3}T_{p^{e-3}})( \Sh_t(f))\\
&= (T_p T_{p^e} -\chi(p^2)( p^{k-3}T_p T_{p^{e-2}} + p^{k-2}T_{p^{e-1}})+ \chi(p^4)p^{2k-5}T_{p^{e-3}})(\Sh_t f)\\
&= (T_{p^{e+1}} - \chi(p^2)p^{k-3}(T_{p^{e-1}} + \chi(p^2)p^{k-2}T_{p^{e-3}})+ \chi(p^4)p^{2k-5}T_{p^{e-3}})(\Sh_t f)\\
&= (T_{p^{e+1}} - \chi(p^2)p^{k-3}T_{p^{e-1}})( \Sh_t f).
\end{split}
\end{equation*}
The first equality uses Theorem~\ref{thm:recur}, third equality follows by using inductive hypothesis 
for $\ell=e$ and $\ell=e-1$, the others follow by using Equation~\eqref{eqn:Heckeprop}.
\end{proof}

We also prove the following proposition, independently of 
the proof of Theorem~\ref{thm:recur}. 
\begin{prop}\label{prop:recur}
Let $p\nmid N$ be a prime and $\ell$ be a positive integer. 
For positive integers $r$ such that 
$1 \le r \le \lfloor\frac{\ell}{2}\rfloor$ we give the following 
recursive construction 
of sequences $A_{r,\ell}(m)$ and $B_{r,\ell}(m)$:
\begin{equation*}
\begin{split}
A_{1,\ell}(m)&= 1, \qquad \qquad A_{r,\ell}(m)= A_{r-1,\ell}(m) - {{\ell-2(r-1)}\choose{m-(r-1)}}A_{r-1,\ell}(r-1); \\
%A_2(m)&= A_1(m) - {{\ell-2}\choose{m-1}}A_1(1) \qquad B_2(m) = B_1(m) - {{\ell-2}\choose{m-1}}B_1(1)\\
%A_3(m)&= A_2(m) - {{\ell-4}\choose{m-2}}A_2(2) \qquad B_3(m) = B_2(m) - {{\ell-4}\choose{m-2}}B_2(2)\\
%\vdots\\
B_{1,\ell}(m)&= {{\ell}\choose{m}}-1, \qquad B_{r,\ell}(m)= B_{r-1,\ell}(m) - {{\ell-2(r-1)}\choose{m-(r-1)}}B_{r-1,\ell}(r-1). 
\end{split}
\end{equation*}
Let $\alpha_{r,\ell} = A_{r,\ell}(r)$ and $\beta_{r,\ell} = B_{r,\ell}(r)$. 
Then 
\[
 T_{p^{2\ell}} = (T_{p^2})^\ell - \sum_{r=1}^{\lfloor\frac{\ell}{2}\rfloor} \chi(p^{2r})
(\alpha_{r,\ell} p^{r(k-2)-1} + \beta_{r,\ell} p^{r(k-2)})(T_{p^2})^{\ell-2r}
\]
as Hecke operators in $\TT_{k/2}$.
\end{prop}
\begin{proof}
Let $f = \sum_{n=0}^{\infty} a(n)q^n \in M_{k/2}(N,\chi)$. Our strategy will be to 
compare the $n$-th coefficient of action of the above operators on $f$ on both sides.
Substituting the $q$-expansion of $f$ in 
Equation~\eqref{eqn:A} and using Lemma~\ref{lem:Gauss} on Gauss sums we obtain 
\[T_{p^{2\ell}} f =I_0+ I_{2\ell} + \sum_{\substack{\nu=1\\ \nu \text{odd}}}^{2\ell-1} I_{\nu}^{\text{odd}} +
\sum_{\substack{\nu=1\\ \nu \text{even}}}^{2\ell-1} I_{\nu}^{\text{even}} \] 
where
\begin{equation*}
\begin{split}
I_0 &=\chi(p^{2\ell})p^{(k-2)\ell}\sum_{n=0}^{\infty} a(n/p^{2\ell}) q^n, \qquad I_{2\ell}= \sum_{n=0}^{\infty} a(np^{2\ell}) q^n \\
I_{\nu}^{\text{odd}} &= \chi(p^{2\ell}-\nu)p^{(\frac{k}{2}-1)(2\ell-\nu)-\frac{1}{2}} \epsilon_p^{k+1}\kro{-1}{p}
\sum_{\substack{n=0 \\ p^{2\ell-\nu-1}\mid n}}^{\infty} a(n/p^{2\ell-2\nu})\kro{n/p^{2\ell-\nu-1}}{p} q^n\\
I_{\nu}^{\text{even}} &= \chi(p^{2\ell}-\nu)p^{(\frac{k}{2}-1)(2\ell-\nu)-1}(
\sum_{\substack{n=0 \\ p^{2\ell-\nu}\mid n}}^{\infty} a(n/p^{2\ell-2\nu})(p-1)q^n-
\sum_{\substack{n=0 \\ p^{2\ell-\nu-1}\| n}}^{\infty} a(n/p^{2\ell-2\nu})q^n).
\end{split}
\end{equation*}
Let $n$ be a positive integer with $p^{2(\ell-1)} \mid n$.
We can write the $n$-th coefficient of $T_{p^{2}}^\ell f$ as 
\begin{equation*}
\begin{split}
& a(np^{2\ell}) + \sum_{m=1}^{\ell-1} {{\ell}\choose{m}} \chi(p^{2m})p^{(k-2)m}a(np^{2\ell-4m}) + \\
& \chi(p^{2\ell-1}) {\kro{-1}{p}}^\frac{k-1}{2}\kro{n/p^{2\ell-2}}{p}p^{\frac{k-3}{2} +(k-2)(\ell-1)}a(n/p^{2\ell-2}) + 
\chi(p^{2\ell})p^{(k-2)\ell}a(n/p^{2\ell}). 
\end{split}
\end{equation*}
Thus the $n$-th coefficient of $T_{p^{2}}^\ell f-T_{p^{2\ell}} f$ is 
\[\sum_{m=1}^{\ell-1} \left({{\ell}\choose{m}}-1\right) \chi(p^{2m})p^{(k-2)m}a(np^{2\ell-4m}) +
\sum_{m=1}^{\ell-1}\chi(p^{2m})p^{(k-2)m-1}a(np^{2\ell-4m}).\]
We want to subtract a suitable multiple of $T_{p^{2}}^{\ell-2} f$ from the above so as to remove 
the terms involving $a(np^{2\ell-4})$ and $a(np^{4-2\ell})$, thereby reducing the number of terms 
in the above sum. Indeed we obtain that the $n$-th coefficient of \\
$(T_{p^{2}}^\ell -T_{p^{2\ell}} - \chi(p^2)(p^{k-3} + (\ell-1)p^{k-2})T_{p^{2}}^{\ell-2} )f$ is 
\begin{equation*}
\begin{split}
& \sum_{m=2}^{\ell-2} \left(1- {{\ell-2}\choose{m-1}}\right) \chi(p^{2m})p^{(k-2)m-1}a(np^{2\ell-4m}) + \\
&\sum_{m=2}^{\ell-2} \left({{\ell}\choose{m}}-1-(\ell-1){{\ell-2}\choose{m-1}}\right)
\chi(p^{2m})p^{(k-2)m}a(np^{2\ell-4m}).
\end{split}
\end{equation*}
We iterate this process of subtracting suitable multiples of 
$T_{p^{2}}^{\ell-2r} f$ which leads us 
to the recursive formulae for $\alpha_{r,\ell}$ and $\beta_{r,\ell}$.
\end{proof}

We obtain the following combinatorial result as a corollary of Theorem~\ref{thm:recur} and 
Proposition~\ref{prop:recur}
\begin{cor}
Keeping the notation as in the previous proposition 
we get the following combinatorial identities for 
$2 \le r \le \lfloor\frac{\ell}{2}\rfloor-1$:
\[\alpha_{r-1,\ell-2} + \alpha_{r,\ell} - \alpha_{r,\ell-1} =0,
\qquad
\beta_{r-1,\ell-2} + \beta_{r,\ell} - \beta_{r,\ell-1} =0.\]
\end{cor}
\begin{proof}
Let $p \nmid N$ be any prime. We substitute the formula for $T_{p^{2\ell}}$ given by Proposition~\ref{prop:recur} 
in the identity of Theorem~\ref{thm:recur}, 
\[
T_{p^{2\ell+2}} - T_{p^2} T_{p^{2\ell}} + \chi(p^2)p^{k-2}T_{p^{2\ell-2}} = 0
\]
to obtain
\begin{equation*}
\begin{split}
& - \sum_{r=2}^{\lfloor\frac{\ell}{2}\rfloor} \chi(p^{2r})
(\alpha_{r,\ell} p^{r(k-2)-1} + \beta_{r,\ell} p^{r(k-2)})(T_{p^2})^{\ell-2r} \\
& + \sum_{r=2}^{\lfloor\frac{\ell-1}{2}\rfloor} \chi(p^{2r})
(\alpha_{r,\ell-1} p^{r(k-2)-1} + \beta_{r,\ell-1} p^{r(k-2)})(T_{p^2})^{\ell-2r}\\
& - \sum_{r=2}^{\lfloor\frac{\ell-2}{2}\rfloor+1} \chi(p^{2r})
(\alpha_{r-1,\ell-2} p^{r(k-2)-1} + \beta_{r-1,\ell-2} p^{r(k-2)})(T_{p^2})^{\ell-2r} = 0.
\end{split}
\end{equation*}
It is clear, with fixed $\ell$ and varying $r$, that the operators 
$(T_{p^2})^{\ell-2r}$ are linearly independent elements 
of $\TT_{k/2}$ and hence 
\[
-\alpha_{r,\ell} + \alpha_{r,\ell-1} -\alpha_{r-1,\ell-2} + 
(\beta_{r,\ell} + \beta_{r,\ell-1} -\beta_{r-1,\ell-2})p=0.
\]
Since this holds for any prime $p$ with $p \nmid N$ the above corollary follows.

\end{proof}

\section{Generators for the Hecke Action}

\begin{thm}\label{thm:heckgen}
Let $k$, $N$ be positive integers with $k\geq 3$ odd,
and $4 \mid N$. Let $\chi$ be a Dirichlet character
modulo $N$. Let $N^\prime=N/2$.
Let $\TT$ be the restriction of Hecke algebra $\TT_{k-1}$ to 
$S_{k-1}(N^\prime,\chi^2)$ and suppose 
$\TT$ is generated as a $\Z$-module
by the Hecke operators $T_{i}$ for $i \le r$.
Then the Hecke operators $T_{i^2}$ for $i \leq r$
generate $\TT_{k/2}^\perp$ as a $\Z[\zeta_{\varphi(N)}]$-module.
%the restriction of Hecke algebra $\TT_{k/2}$ to 
%$S_{k/2}^\perp(N,\chi)$ as a $\Z[\zeta_N]$-module. 
In particular, $f \in S_{k/2}^\perp(N,\chi)$
is an eigenform for all Hecke operators if and only if it is
an eigenform for $T_{i^2}$ for $i \le r$.
\end{thm}

\begin{proof}
Let $n$ be a positive integer with prime factorization $n= p_1^{n_1}p_2^{n_2}\cdots p_s^{n_s}$. 
Let $f \in S_{k/2}^\perp(N,\chi)$. Let $t$ be a square-free positive integer. Using 
Theorem~\ref{thm:recur} or Proposition~\ref{prop:recur}, for any prime $p$ and a positive 
integer $\ell$ we can express the action of $T_{{p}^{2\ell}}$ as
\begin{equation}\label{eqn:Tpl}
T_{{p}^{2\ell}} = \sum_{j=0}^{\ell} \gamma_{j}T_{p^2}^{j}, \qquad \gamma_{j} \in \Z[\zeta_{\varphi(N)}].
\end{equation}
Note that in the above expression $\gamma_\ell=1$ and hence the Hecke operators 
$T_{p^2}^{j}$ with $ 1 \le j \le \ell$ generates the same $\Z[\zeta_{\varphi(N)}]$-module as does the 
Hecke operators $T_{{p}^{2j}}$ with $ 1 \le j \le \ell$.
Thus we have 
\begin{equation}\label{eqn:del1}
\begin{split}
\Sh_t(T_{n^2}f) &= \Sh_t( T_{p_1^{2n_1}}T_{p_2^{2n_2}}\cdots T_{p_s^{2n_s}} f)\\
&= \Sh_t\left(\left(\sum_{j_1=0}^{n_1} \gamma_{j_1}T_{p_1^2}^{j_1}\right) \cdots 
\left(\sum_{j_s=0}^{n_s} \gamma_{j_s}T_{p_s^2}^{j_s}\right)f\right)\\
&= \left(\sum_{j_1=0}^{n_1} \gamma_{j_1} T_{p_1}^{j_1}\right) \cdots 
\left(\sum_{j_s=0}^{n_s} \gamma_{j_s} T_{p_s}^{j_s}\right)\left(\Sh_t f \right)\\
&= \sum_{i=1}^{r} \delta_i T_i (\Sh_t f), 
\end{split}
\end{equation}
where the last equality follows since the $T_i$, with $1\le i\le r$, generate $\TT$ as a $\Z$-module, while 
the second last equality follows by Proposition~\ref{prop:commute}.

Recall from Proposition~\ref{prop:Heckeprop1}, for any prime $q$ and a positive integer $\ell$, 
the action of Hecke operator $T_{q^\ell}$ on $S_{k-1}(N^\prime,\chi^2)$ can be expressed as 
\[T_{q^\ell} = \sum_{j=0}^{\ell}\alpha_j T_q^j, \qquad \alpha_j \in \Z[\zeta_{\varphi(N^\prime)}] \subset \Z[\zeta_{\varphi(N)}]. \]

Let $1\le i \le r$ has prime factorization $i = q_1^{m_1}q_2^{m_2}\cdots q_v^{m_v}$. 
Then each term $T_i (\Sh_t f)$ in Equation~\eqref{eqn:del1} can be written as 
\begin{equation}\label{eqn:del2}
\begin{split}
 T_i (\Sh_t f) &= T_{q_1^{m_1}}T_{q_2^{m_2}}\cdots T_{q_v^{m_v}}(\Sh_t f)\\
&=  \left(\sum_{j_1=0}^{m_1} \alpha_{j_1}T_{q_1}^{j_1}\right) \cdots  
\left(\sum_{j_v=0}^{m_v} \alpha_{j_v}T_{q_v}^{j_v}\right) (\Sh_t f)\\
%&=  \left(\left(\sum_{j_{i_1}=0}^{m_{i_1}} \alpha_{j_{i_1}}(T_{q_{i_1}})^{j_{i_1}}\right) \cdots  
%\Sh_t \left(\sum_{j_{i_v}=0}^{m_{i_v}} \alpha_{j_{i_v}}(T_{q_{i_v}^2})^{j_{i_v}}\right)f\right)\\
&=  \Sh_t \left(\left(\sum_{j_1=0}^{m_1} \alpha_{j_1}T_{q_1^2}^{j_1}\right) \cdots  
\left(\sum_{j_v=0}^{m_v} \alpha_{j_v} T_{q_v^2}^{j_v}\right)f\right)\\
&=  \Sh_t \left(\left(\sum_{j_1=0}^{m_1} \beta_{j_1} T_{q_1^{2j_1}}\right) \cdots  
\left(\sum_{j_v=0}^{m_v} \beta_{j_v} T_{q_v^{2j_v}}\right)f\right)\\
&= \Sh_t \left( \sum_{j=1}^{i} A_j T_{j^2} f \right),
\end{split}
\end{equation}
where $A_j \in \Z[\zeta_{\varphi(N)}]$. In the above equalities we repeatedly use Proposition~\ref{prop:commute} 
and Equation~\eqref{eqn:Tpl}. For the second last equality we use the remark below Equation~\eqref{eqn:Tpl}. 
Now using Equations \eqref{eqn:del1} and \eqref{eqn:del2} we get 
\[\Sh_t(T_{n^2} f) = \Sh_t \left( \sum_{i=1}^{r} B_i T_{i^2} f \right), \qquad B_i \in \Z[\zeta_{\varphi(N)}].\]
Since this is true for all positive square-free integers $t$, using Lemma~\ref{lem:shker} we deduce that 
\[ T_{n^2} f =  \sum_{i=1}^{r} B_i T_{i^2} f. \]
Hence $T_{i^2}$ with $i \leq r$ generate $\TT_{k/2}^\perp$ as a $\Z[\zeta_{\varphi(N)}]$-module.
%generate the restriction of $\TT_{k/2}$ to $S_{k/2}^\perp(N,\chi)$ as 
%a $\Z[\zeta_N]$-module.
\end{proof}

We shall need the following theorem which is a consequence of 
Sturm's bound \cite{Sturm}.

\begin{thm} (Stein \cite[Theorem 9.23]{Stein}) \label{thm:LS}
Suppose $\Gamma$ is a congruence subgroup that contains
$\Gamma_1(N)$. Let
\[
r= \frac{km}{12}-\frac{m-1}{N}, \qquad m=[\SL_2(\Z) : \Gamma].
\]
Then the Hecke algebra
\[
\TT=\Z[\dots,T_n,\dots] \subset \End(S_k(\Gamma))
\]
is generated as a $\Z$-module by the Hecke operators $T_n$ for
$n \leq r$.
\end{thm}
%Note that $\Gamma_0(N)$ contains $\Gamma_1(N)$,
%so the above theorem holds for $\Gamma=\Gamma_0(N)$.

From Theorem~\ref{thm:LS} we deduce the following corollary.
\begin{cor} 
Let $k$, $N$ be positive integers with $k\geq 3$ odd,
and $4 \mid N$. Let $\chi$ be a Dirichlet character
modulo $N$. Let $N^\prime=N/2$. Write
\[
m= {N^\prime}^2 \prod_{p \mid N^\prime} 
\left( 1- \frac{1}{p^2} \right),
\qquad 
R= \frac{(k-1)m}{12}-\frac{m-1}{N^\prime}.
\]
Then $T_{i^2}$ for $i \leq R$
generate $\TT_{k/2}^\perp$ as a $\Z[\zeta_{\varphi(N)}]$-module. 
In particular the set of operators $T_{p^2}$ for primes $p \leq R$ forms a 
generating set as an algebra.
%Moreover, $f \in S_{k/2}(N,\chi)$
%is an eigenform for all Hecke operators if and only if it is
%an eigenform for $T_{p^2}$ for $p\leq R$.
%\end{cor}
%\begin{cor}\label{cor:sturm}
%With the same hypothesis as in the above corollary, further suppose that 
Moreover, if 
$\chi$ is a quadratic character, then the same result holds as above with 
\[
m = N^\prime \prod_{p \mid N^\prime} \left( 1+ \frac{1}{p} \right),
\qquad 
R= \frac{(k-1)m}{12}-\frac{m-1}{N^\prime}.
\]
\end{cor}
\begin{proof}
Note that $S_{k-1}(N^\prime, \chi^2) \subset S_{k-1}(\Gamma_1(N^\prime))$. The first part of the corollary 
follows by applying Theorem~\ref{thm:heckgen} and Theorem~\ref{thm:LS} to the 
congruence subgroup $\Gamma_1(N^\prime)$ and using the formula for $[\SL_2(\Z):\Gamma_1(N^\prime)]$ that 
can be found for example in \cite[Page 14]{DS}.

Now suppose $\chi$ is a quadratic character. Then $S_{k-1}(N^\prime,\chi^2)=S_{k-1}(N^\prime)$. So we apply 
Theorem~\ref{thm:LS} to the group $\Gamma_0(N^\prime)$ and we now use the formula for 
$[\SL_2(\Z):\Gamma_0(N^\prime)]$.
\end{proof}

\end{document}